\documentclass[oneside,english,11pt]{amsart}
\pdfoutput=1
\usepackage[T1]{fontenc}
\usepackage[latin9]{inputenc}
\usepackage[letterpaper]{geometry}
\usepackage{amsthm}
\usepackage{amssymb}

\numberwithin{equation}{section} 
\numberwithin{figure}{section} 
\theoremstyle{plain}
\theoremstyle{plain}
\newtheorem{thm}{Theorem}
  \theoremstyle{plain}
  \newtheorem{prop}[thm]{Proposition}
  \theoremstyle{definition}
  \newtheorem{defn}[thm]{Definition}
  \theoremstyle{plain}
  \newtheorem{lem}[thm]{Lemma}

\usepackage{babel}

\begin{document}

\title{Vertices of degree $k$ in edge-minimal, $k$-edge-connected graphs}

\author{Carl Kingsford}

\address{C. Kingsford, Department of Computer Science and Institute for Advanced
Computer Studies, University of Maryland, College Park, MD, USA}

\email{carlk@cs.umd.edu}

\thanks{C.K. was partially supported by NSF grant IIS-0812111.}

\author{Guillaume Mar\c{c}ais}

\address{G. Mar\c{c}ais, Program in Applied Mathematics \& Statistics and Scientific
Computation, University of Maryland, College Park, MD, USA}

\email{guillaume@marcais.net}

\date{\today}

\maketitle
Halin~\cite{Halin1969} showed that every edge-minimal, $k$-\textit{vertex}-connected
graph has a vertex of degree $k$. In this note, we prove the analogue
to Halin's theorem for edge-minimal, $k$-\textit{edge}-connected
graphs:
\begin{thm}
\label{thm:Let--be}Let $G$ be an edge-minimal, $k$-edge-connected
graph. Then there are two nodes of degree $k$ in $G$.
\end{thm}
To prove Theorem~\ref{thm:Let--be}, we first establish a link between
edge-minimal, $k$-edge-connected graphs and exactly $k$-edge-connected
graph \cite{KingsfordMarcais:2009} (Definition~\ref{def:Exact-connectivity}
and Proposition~\ref{pro:G^k-exactly-connected}). The theorem is
proved in the case of $G$ being an exactly $k$-edge-connected graph
(Proposition~\ref{pro:Exact-connected-2-deg-k}) and then transfered
to edge-minimal graphs. Throughout, all graphs considered are multigraphs
and all sets are multisets.

\subsection*{The edge-connectivity equivalence relation. }
\begin{prop}
The edge version of Menger's theorem holds for multigraphs.\end{prop}
\begin{proof}
Any multi-edge in the graph $G$ can be replaced by a length $2$
path to get a graph $G'$. There is then an obvious bijection between
the paths in the original multigraph $G$ and the resulting graph
$G'$. The result of the Menger theorem on $G'$ can be immediatly
applied to $G$.\end{proof}
\begin{defn}
Let $v_{1}R^{k}v_{2}$ be the relation between nodes $v_{1}$ and
$v_{2}$ that holds if $v_{1}=v_{2}$ or there are $k$ edge-disjoint
paths between $v_{1}$ and $v_{2}$.\end{defn}
\begin{prop}
$R^{k}$ is an equivalence relation.\end{prop}
\begin{proof}
$R^{k}$ is by definition reflexive and symmetric. Let $G$ be a graph,
let $v_{1},v_{2},v_{3}\in G$ satisfy $v_{1}R^{k}v_{2}$ and $v_{2}R^{k}v_{3}$.
If $v_{1}=v_{2}$, $v_{2}=v_{3}$ or $v_{1}=v_{3}$, then transitivity
is obvious. Suppose that $v_{1}$, $v_{2}$ and $v_{3}$ are distinct
vertices and let $S$ be an edge set of cardinality $k-1$. There
are $k$ edge-disjoint paths from $v_{1}$ to $v_{2}$ in $G$ so,
$v_{1}$ and $v_{2}$ are still connected in $G\setminus S$ and so
are $v_{2}$ and $v_{3}$. So we have a path $v_{1}\rightsquigarrow v_{2}\rightsquigarrow v_{3}$
in $G\setminus S$. By Menger's theorem, the set of paths between
$v_{1}$ and $v_{3}$ in $G$ has cardinality at least $k$ and $v_{1}R^{k}v_{3}$,
and $R^{k}$ is transitive.\end{proof}
\begin{prop}
\label{pro:k-1 paths in R^k}There are at most $k-1$ edge-disjoint
paths between two equivalency classes of $R^{k}$.\end{prop}
\begin{proof}
Suppose we have $k$ edge-disjoint paths, $p_{i}$, $1\le i\le k$,
between $C_{1}$ and $C_{2}$, two distinct equivalency classes of
$R^{k}$ in graph $G$. Let $v_{i}^{1}$, $1\le i\le k$ and $v_{i}^{2}$,
$1\le i\le k$ be the endpoints in $C_{1}$ and $C_{2}$ respectevily
of $k$ edge-disjoint paths. (Note that the $v_{i}^{1}$ are not all
necessary distinct. This is true for the $v_{i}^{2}$ as well.) Let
$S$ be a set of $k-1$ edges. Then, in $G\setminus S$, at least
one path, say $p_{i_{0}}$, was not disconnected. Because $v_{1}^{1}R^{k}v_{i_{0}}^{1}$,
$v_{1}^{1}$ and $v_{i_{0}}^{1}$ are not disconnected. Similarly,
$v_{1}^{2}$ and $v_{i_{0}}^{2}$ are not disconnected. So we have
a path $v_{1}^{1}\rightsquigarrow v_{i_{0}}^{1}\rightsquigarrow v_{i_{0}}^{2}\rightsquigarrow v_{1}^{2}$.
By Menger's theorem, $v_{1}^{1}R^{k}v_{1}^{2}$, which is a contradiction.
\end{proof}

\subsection*{Relationship between exact connectivity and minimality.}
\begin{prop}
\label{pro:edge-minimal-adjacency}Let $G=(V,E)$ be a $k$-edge-connected
graph. $G$ is edge minimal if and only if for any adjacent vertices
$(u,v)\in E$, there are at most $k$ edge disjoint $u-v$ paths.\end{prop}
\begin{proof}
Let $e=(u,v)\in E$ such that there are $>k$ edge-disjoint $u-v$
paths. In $G-e$, there are at least $k$ edge-disjoint $u-v$ paths.
Let $x$ and $y$ be any two vertices, not necessarily distinct from
$u$ and $v$. There are $k$ edge-disjoint $x-u$ paths in $G$.
So, depending on whether or not edge $e$ is on one of these $x-u$
paths, in $G-e$ there are either $k$ edge-disjoint $x-u$ paths,
or there are $k-1$ edge-disjoint $x-u$ paths and a $x-v$ path edge
disjoint from the $k-1$ $x-u$ paths. There is a similar situation
in $G-e$ between $y$, $v$ and $u$. Let $S$ be a set of $k-1$
distinct edges of $G-e$. In $G-e-S$, there is a $x-u$ or $x-v$
path, a $y-u$ or $y-v$ path and a $u-v$ path. Hence there is a
$x-y$ path and $S$ is not a separating set. By Menger's theorem,
$G-e$ is $k$-edge-connected, and $G$ is not edge minimal. Conversely,
suppose that for any edge $e=(u,v)\in E$ there are at most $k$ edge-disjoint
$u-v$ paths. Then in $G-e$ there are at most $k-1$ edge disjoint
$u-v$ paths. So $G$ is edge minimal.\end{proof}
\begin{defn}
\label{def:Exact-connectivity}A graph $G$ is called \textbf{exactly
$k$-edge-connected} if there are exactly $k$ edge disjoint paths
between any two nodes $u,v\in G$.
\end{defn}

\begin{defn}
Let $G$ be a $k$-edge-connected graph. Define $G^{k}$ to be a graph
where the vertices are the equivalency classes of $R^{k+1}$ on $G$
and there is an edge $(v_{1},v_{2})\in G^{k}$ for every edge $(u_{1},u_{2})\in G$
with $u_{1}\in v_{1}$ and $u_{2}\in v_{2}$.\end{defn}
\begin{prop}
\label{pro:G^k-exactly-connected}If $G$ is an edge-minimal, $k$-edge-connected
graph then $G^{k}$ is not trivial and is exactly $k$-edge-connected.\end{prop}
\begin{proof}
If $G$ is edge minimal, it is not $(k+1)$-edge-connected, so $R^{k+1}$
has more than one equivalency class, and $G^{k}$ is not trivial.
$G^{k}$ is $k$-edge-connected like $G$. By Proposition~\ref{pro:k-1 paths in R^k},
it is exactly $k$-edge-connected.\end{proof}
\begin{prop}
\label{pro:minimal-degree-classes}Let $G$ be an edge-minimal, $k$-edge-connected
graph, and let $v\in G^{k}$ be an equivalence class of $R^{k+1}$.
Then, for any $u\in v$, $\deg_{G^{k}}(v)\ge\deg_{G}(u)$.\end{prop}
\begin{proof}
Let $v\in G^{k}$. By Proposition \ref{pro:edge-minimal-adjacency},
if $u_{1},u_{2}\in v$, then $(u_{1},u_{2})$ is not an edge in $G$.
So every neighbor of $u_{1}$ is not in $v$ and by construction $\deg_{G^{k}}(v)\ge\deg_{G}(u_{1})$.
\end{proof}

\subsection*{Proof of Theorem~\ref{thm:Let--be}.}
\begin{defn}
An edge cut $S$ of a graph $G$ is called \textbf{trivial} if one
of the components of $G\setminus S$ is the trivial graph.
\end{defn}

\begin{defn}
A \textbf{$k$-regular} graph is a graph where all vertices have the
same degree $k$. A \textbf{quasi $k$-regular} graph is a graph where
at most one vertex has a degree different than $k$.\end{defn}
\begin{lem}
An exactly $k$-edge connected graph $G$ which has only trivial cuts
is quasi $k$-regular.\end{lem}
\begin{proof}
Suppose there exists two vertices $u$ and $v$ of degree greater
than $k$. There exists a minmum cut separating $u$ and $v$ and
this cut cannot be trivial.\end{proof}
\begin{defn}
[Vertex splitting]Let $G=(V,E)$ be an exactly $k$-edge connected
graph and $S=\langle V_{1},V_{2}\rangle$ be a non-trivial minimum
cut. Construct $G_{1}=(V_{_{1}}\cup\left\{ x_{1}\right\} ,E_{1})$
and $G_{2}=(V_{2}\cup\left\{ x_{2}\right\} ,E_{2})$ by adding two
new vertices $x_{1}$ and $x_{2}$ attached respectively to $G_{1}$
and $G_{2}$ by $k$ new edges to the vertices adjacent to $S$. Formally,
let $S=\left\{ (u_{i},v_{i})\in V_{1}\times V_{2}:1\le i\le k\right\} $
and\begin{eqnarray*}
E_{1} & = & \left(V_{1}\times V_{1}\cap E\right)\cup\left\{ (x_{1},u_{i}):1\le i\le k\right\} \\
E_{2} & = & \left(V_{2}\times V_{2}\cap E\right)\cup\left\{ (x_{2},v_{i}):1\le i\le k\right\} \,.\end{eqnarray*}
The pair $G_{1},G_{2}$ is called a \textbf{vertex splitting} of $G$
with respect to $S$.\end{defn}
\begin{prop}
\label{pro:Vertex-splitting}Let $G$ be an exactly $k$-edge-connected
graph, and let $S$ be a non-trivial minimum cut. $G_{1}$ and $G_{2}$
obtained from vertex splitting $G$ with respect to $S$ are exactly
$k$-edge-connected.
\end{prop}

\begin{prop}
\label{pro:Exact-connected-2-deg-k}Let $G$ be an exactly $k$-edge-connected
graph. Then there are two vertices of degree $k$ in $G$.\end{prop}
\begin{proof}
We proceed by induction on the number of non-trivial minimum cuts
in $G$. If $G$ has no non-trivial minimum cuts, then it is quasi
$k$-regular and has at least 2 nodes of degree $k$. Let $S=\langle V_{1},V_{2}\rangle$
be a non-trivial minimum cut and let $G_{1}$ and $G_{2}$ be the
vertex splitting graphs induced by $S$. Call $x_{1}$ and $x_{2}$
the new vertices ($V(G_{i})=V_{i}\cup\{x_{i}\}$). Suppose $T$ were
a non-trivial minimum cut in $G_{1}$. Construct a corresponding non-trivial
minimum cut $T'$ in $G$ by changing any edge used by $T$ that is
adjacent to $x_{1}$ to the corresponding edge in $S$. So to any
non-trivial minimum cut of $G_{1}$ or $G_{2}$ corresponds a distinct
non-trivial minimum cut in $G$. But no non-trivial cut in $G_{1}$
or $G_{2}$ corresponds to the cut $S$ (it would be a trivial cut
in $G_{1}$ and $G_{2}$). So both $G_{1}$ and $G_{2}$ have fewer
non-trivial minimum cuts than $G$. By induction $G_{1}$ and $G_{2}$
have 2 nodes of degree $k$, including $x_{1}$ and $x_{2}$. So,
$G$ has the same vertices as $G_{1}$ and $G_{2}$, except for $x_{1}$
and $x_{2}$, with the same degree. Hence it has 2 nodes of degree
$k$.
\end{proof}
\bigskip{}

\begin{proof}
[Proof of Theorem~\ref{thm:Let--be}]By Proposition~\ref{pro:Exact-connected-2-deg-k},
$G^{k}$ has two vertices of degree $k$. Let $v$ be an equivalence
class of $R^{k+1}$. If $v$ contains two distinct vertices $u_{1}$
and $u_{2}$ of $G$, then there are at least $k+1$ edge disjoint
$u_{1}-u_{2}$ paths in $G$ and both $u_{1}$ and $u_{2}$ have a
degree $\ge k+1$. By Proposition \ref{pro:minimal-degree-classes},
$v$ also has degree $\ge k+1$ in $G^{k}$. Hence, the vertices of
degree $k$ in $G^{k}$ correspond to equivalence classes that must
each contain at most one vertex of $G$ of degree $k$. Therefore,
$G$ has two vertices of degree $k$. 
\end{proof}
\bibliographystyle{hplain}
\bibliography{bib}

\bigskip{}

\end{document}